\documentclass[12pt]{article}

\makeatletter
\renewcommand{\@oddfoot}{\hfill \thepage}
\makeatother

\usepackage[T2A]{fontenc}
\usepackage[cp1251]{inputenc}
\usepackage[russian, english]{babel}
\usepackage{mathrsfs}
\usepackage{amssymb}
\usepackage{cite}
\usepackage{amsmath,amsthm,amsfonts}
\usepackage{epstopdf}
\usepackage[dvipdf]{graphicx}
\usepackage{fancyhdr}

\begin{document}

\begin{center}
{\bf SLOW DIFFUSION BY MARKOV RANDOM FLIGHTS}
\end{center}

\begin{center}
Alexander D. KOLESNIK\\
Institute of Mathematics and Computer Science\\
Academy Street 5, Kishinev 2028, Moldova\\
E-Mail: kolesnik@math.md
\end{center}

\vskip 0.2cm

\begin{abstract}
We present a conception of the slow diffusion processes in the Euclidean spaces $\Bbb R^m, \; m\ge 1$, 
based on the theory of random flights with small constant speed that are driven by a homogeneous Poisson process of small rate. 
The slow diffusion conditions that, on long time intervals, lead to the stationary distributions, are given. 
The stationary distributions of slow diffusion processes in some Euclidean spaces of low dimensions, are presented.  
\end{abstract}

\vskip 0.1cm


\vskip 0.2cm


\section{Preliminaries}
\numberwithin{equation}{section}

Various diffusion processes often arise in the nature as well as in many fields of science and technology. These processes are  characterized by the presence of some source of substance (or heat, mass, energy, electric charge, etc.) concentrated at an initial point or in a compact set and spreading the substance outwards. The attempts of describing such processes have very long history. The first one was, apparently, undertaken by W. Thomson (Lord Kelvin) \cite{lord} who, still in the middle of the 19-th century, has suggested to describe the propagation of electric signals in the transatlantic cable by the damped wave (telegraph) equation. A similar telegraph equation has been used by V.A. Fock \cite{fock} for describing the diffusion of a light ray in a homogeneous medium. The revolutionary idea to interpret a diffusion process as the substance transport by some chaotically moving Brownian particles, was almost simultaneously put forward by A. Einstein \cite{ein}, M. Smoluchowski \cite{smol} and developed then by G.I. Taylor \cite{taylor}. The underlying idea of transport by Brownian particles leads to the one-dimensional parabolic heat equation (see equation (\ref{heat}) below), whose fundamental solution (i.e. the Green's function) represents the density of a Brownian particle's position on the line, that yields the respective distribution of the substance mass. This one-dimensional Einstein-Smoluchowski's model can easily be generalized for the case of transport in the Euclidean space of arbitrary dimension. Namely, in this case the one-dimensional heat equation should be replaced by its multidimensional counterpart whose fundamental solution yields the density of a Brownian particle moving in the space. 

After the appearance of the Einstein-Smoluchowski's model it was extensively applied to describe various real phenomena in many fields of science and technology. It was noted that the theoretical calculations based on this model agree well with experimental data if the speed of the process is sufficiently large. However, it was fairly quickly discovered that in the cases when the speed is small, this model yields unsatisfactory results. This fact can easily be explained if we take into account that the Einstein-Smoluchowski's model is based on the Brownian motion, that has the infinite speed of propagation and the Brownian particle changes its direction infinitely many times per unit of time. These fairly artificial and somewhat paradoxical properties of the Brownian motion are the main reason of the unsatisfactoriness of the classical Einstein-Smoluchowski's infinite-velocity transport model for describing the processes with slow speed of propagation. 

That is why an alternative approach for describing the transport processes with a finite speed of propagation, was suggested. In contrast to the Einstein-Smoluchowski's model, it was assumed that in this case the transport is carried out by the particles that move at {\it finite} speed and change their directions a {\it finite} number of times. A random walk in a plane with steps of a fixed finite length and random directions was first introduced and studied in the works by K. Pearson \cite{pearson1,pearson2}. The three-dimensional counterpart of such random walk was presented by J.W.S. Rayleigh  \cite{ray} (for more details on these stochastic motions see also \cite[Chapter 4]{pin}). It is clear that, since the steps have fixed length, the changes of direction occur at discrete (non-random) time instants. These works are the milestones that contain a hidden hint of the finite-velocity transport by the particles whose motion represents a stochastic process with finite speed of propagation. This idea was realized in the pioneer works by S. Goldstein \cite{gold} and M. Kac \cite{kac} who were the first to associate the finite-velocity random walk and the telegraph equation. The subject of their interest is the stochastic motion of a particle that, at the initial time $t=0$, starts from the origin $x=0$ of the real line $\Bbb R^1$ and moves with a constant finite speed $c$ by alternating two possible directions (positive and negative) at random time instants that form a homogeneous Poisson flow of rate $\lambda>0$. Since the speed $c$ is finite then, at arbitrary time $t>0$, the support of the distribution of the particle's position $X(t)$ represents the closed interval $[-ct, ct]$. In these works it was shown that the probability density $p(x,t), \; x\in\Bbb R^1, t>0,$ of the process $X(t)$ at arbitrary time $t>0$ satisfies the following hyperbolic partial differential equation of second order with constant coefficients: 
\begin{equation}\label{tel}
\frac{\partial^2 p(x,t)}{\partial t^2} + 2\lambda \; \frac{\partial p(x,t)}{\partial t} = c^2 \; 
\frac{\partial^2 p(x,t)}{\partial x^2} ,  
\end{equation}
(which is known as the {\it telegraph} or {\it damped wave} equation) and can be found by solving (\ref{tel}) with the initial conditions: 
\begin{equation}\label{init}
p(x,t)|_{t=0} = \delta(x), \qquad \left. \frac{\partial p(x,t)}{\partial t}\right|_{t=0} = 0 ,  
\end{equation}
where $\delta(x)$ is the Dirac delta-function. From this fact it follows that the density $p(x,t)$ is the fundamental solution (the Green's function) to the telegraph equation (\ref{tel}). A telegraph equation with constant coefficients is one of the classical equation of mathematical physics and its fundamental solution is well known (see, for instance, \cite{morse}). In this particular 
case, the fundamental solution of equation (\ref{tel}) is given by the formula (see \cite[Theorem 2.5]{kolrat}): 

\begin{equation}\label{dens1}
\aligned
p(x,t) & = \frac{e^{-\lambda t}}{2} \left[ \delta(ct+x) +
\delta(ct-x) \right]\\
& + \frac{e^{-\lambda t}}{2c} \left[ \lambda I_0\left(
\frac{\lambda}{c} \sqrt{c^2t^2-x^2}\right) +
\frac{\partial}{\partial t} \; I_0\left( \frac{\lambda}{c}
\sqrt{c^2t^2-x^2}\right)  \right] \Theta(ct-\vert x\vert) , \\
& \qquad\qquad\qquad x\in\Bbb R^1, \qquad t>0, \endaligned
\end{equation}
where
$$I_0(z) = \sum_{k=0}^{\infty} \frac{1}{(k!)^2} \left( \frac{z}{2} \right)^{2k}$$
is the modified Bessel function of order zero and $\Theta(x)$ is the Heaviside unit-step function
\begin{equation}\label{heaviside}
\Theta(x) = \left\{ \aligned & 1, \qquad x>0,\\
& 0, \qquad x\le 0. \endaligned \right. 
\end{equation}
Note that the telegraph equation can also be considered in a more general context as a particular case of the Maxwell equation 
(see \cite[Section 2, subsection 6]{vlad}). The first term of (\ref{dens1}) represents the density of the singular component of the distribution of $X(t)$ (in the sense of generalized functions) concentrated at the two terminal points $\pm ct$ of the interval 
$[-ct, ct]$, while the second term is the density of the absolutely continuous component of the distribution concentrated 
in the open interval $(-ct, ct)$. 

The hyperbolicity of the governing telegraph equation (\ref{tel}) is the main peculiarity of this approach initiated by S. Goldstein \cite{gold} and M. Kac \cite{kac} and developed then in the works by C.R. Cattaneo \cite{catt}, P. Vernott \cite{vern}, E.V. Tolubinsky \cite{tolub} and some other authors. This fact can obviously be explained by the finite velocity of the particle's motion because, from the general PDEs' theory, it is well known that the processes with a finite speed of propagation are described just by hyperbolic partial differential equations. This approach, that leads to the {\it hyperbolic} telegraph equation (\ref{tel}), enables to avoid many strangenesses and paradoxes inherent in the Einstein-Smoluchowski's infinite-velocity model governed by the {\it parabolic} heat equation. These reasonings allow us to assert that the diffusion processes, whose speed of propagation is not too big, can effectively be simulated by the Goldstein-Kac stochastic motion rather than the Einstein-Smoluchowski's model. 

The remarkable fact noted by M. Kac \cite{kac} is that, if both the speed of motion $c$ and the intensity of alternating 
the directions $\lambda$ tend to infinity in such a way that the following scaling condition fulfils
\begin{equation}\label{KacCond}
c\to\infty, \qquad \lambda\to\infty, \qquad \frac{c^2}{\lambda}\to\rho^2, 
\end{equation}
then telegraph equation (\ref{tel}) turns into the Einstein-Smoluchowski's heat equation 
\begin{equation}\label{heat}
\frac{\partial u(x,t)}{\partial t} = \frac{\rho^2}{2} \; \frac{\partial^2 u(x,t)}{\partial x^2} .  
\end{equation}
As it should be expected, under the Kac's scaling condition (\ref{KacCond}), the transition density (\ref{dens1}) transforms into 
the fundamental solution of the parabolic heat equation (\ref{heat}) (see \cite[Section 2.6]{kolrat}), that is, into the transition density of the one-dimensional homogeneous Brownian motion with zero drift and diffusion coefficient $\rho^2$. Therefore, Kac's condition (\ref{KacCond}) can be interpreted as the {\it fast diffusion condition}. 

After the appearance of the Goldstein-Kac telegraph process many efforts were made to generalize it for the case of finite-velocity motions in the high-dimensional Euclidean spaces $\Bbb R^m, \; m\ge 2$. This problem, fisrt formulated by M. Kac, sounds as follows: {\it Can the finite-velocity random motion in higher dimensions be described by the multidimensional telegraph equation similarly to the one-dimensional case?} This question has become the subject of intense discussions among specialists. Some of them, based on the above mentioned analogy between the telegraph process and the one-dimensional Brownian motion, tried to describe such stochastic   motions by means of the multidimensional telegraph equation similar to (\ref{tel}), in which the spatial operator $\partial^2/\partial x^2$ was formally replaced by the Laplace operator $\Delta$ of respective dimension. Other researchers claimed that such a formal replacement was highly doubtful and unjustified from the mathematical point of view. Bartlett \cite[p. 705]{bart} wrote that "such equivalence is more doubtful in the multidimensional case". Tolubinsky \cite[p. 49]{tolub} has characterized such attempts as "unjustified". The final solution of this problem was given in \cite{kolpin} where it was shown that the multidimensional finite-velocity random motions (also called the Markov random flights or, in a more general context, the random evolutions \cite{pin}) are driven by much more complicated equations than the telegraph one, namely, by the hyperparabolic operators represented by the infinite series composed of the integer powers of the telegraph and Laplace operators.  Nevertheless, the Kac's scaling condition (\ref{KacCond}) keeps working in the multidimensional case as well providing the transformation of the hyperparabolic operators into the heat operator of respective dimension (see \cite[Theorem 2]{kolpin}). Moreover, the more strong result states that, under Kac's condition (\ref{KacCond}), the transition density of the symmetric Markov random flight in the Euclidean space $\Bbb R^m$ of arbitrary dimension $m\ge 2$ is convergent to the transition density of the $m$-dimensional homogeneous Brownian motion  with zero drift and diffusion coefficient $2\rho^2/m$ (see \cite[Theorem 4]{kol3}). All these facts lead us to the conclusion that the Kac's scaling condition (\ref{KacCond}) has an universal character in the space of arbitrary dimension and, therefore, it can be treated as the fast diffusion condition indeed. 
Note that no scaling conditions, similar to (\ref{KacCond}), were obtained so far for the finite-velocity stochastic motions driven by various non-Markovian processes studied in a series of recent works \cite{cas,lecaer1,lecaer2,lecaer3,letac}. 

While the fast diffusion processes are of a certain interest in studying many real phenomena, there are numerous important dynamic processes that can be interpreted as the slow diffusion ones. Such processes are characterized first of all by the slow and super-slow speed of propagation and their evolution can last very long time (months, years and even decades). They are of a special importance 
due to their numerous applications in physics, chemistry, biology, environmental science and some other fields (see, for instance,  \cite{aron,chung,giona1,giona2,peters,valdes,wald,zhong} and the bibliographies therein). The mathematical methods applied in some of these works are mostly based on differential equations \cite{aron,zhong}.  

In this article we present a conception of the slow diffusion processes based on the theory of Markov random flights in the Euclidean spaces developed in recent years. The important point of this approach relies on the fact that the main probabilistic characteristics of such processes (transition densities, for instance) in some spaces of low dimensions were obtained in explicit forms (see \cite{kol2,kol5,kolors,mas,sta1,sta2}) that can be used for deriving the stationary distributions of slow diffusion processes.  
The crucial point is to determine the conditions on the parameters of the motion under which the random flight generates a {\it slow} diffusion process (similarly like the Kac's conditions (\ref{KacCond}) generates a {\it fast} diffusion one). However, in contrast to the fast diffusion condition (\ref{KacCond}), we should determine the appropriate conditions not only for the speed of motion and the intensity of switchings, but also for the time variable. The slow speed of propagation implies that we should consider the process on large time intervals and this leads to the respective stationary distributions.  

The article is organized as follows. In Section 2 we give a description of the multidimensional symmetric Markov random flight and present a condition (in fact, a set of conditions) under which it generates a slow diffusion process in the Euclidean space of arbitrary dimension. This means that such slow diffusion condition has an universal character like the Kac's fast diffusion one (\ref{KacCond}). This condition connects the slow speed of propagation and the small rate of Poisson switchings through the time and, on long time intervals, lead to stationary distributions. 
In Section 3 we apply this slow diffusion condition to the known distributions of the symmetric Markov random flights in the Euclidean spaces $\Bbb R^1, \Bbb R^2, \Bbb R^4, \Bbb R^6$ and obtain closed-form expressions for their stationary distributions. In the space $\Bbb R^3$ such a stationary distribution is given in an asymptotic form under the additional condition of the smallness of the mean value (the expectation) of the number of switchings on the respective time interval.

\section{Random Flights and Slow Diffusion Condition}
\numberwithin{equation}{section} 

The multidimensional counterpart of the Goldstein-Kac telegraph process is the symmetric Markov random flight represented by the stochastic motion of a particle that, at the initial time instant $t=0$,  starts from the origin $\bold 0 = (0,\dots,0)$ of the Euclidean space $\Bbb R^m, \; m\ge 2,$ and moves with some constant finite speed $c$. The initial direction is a random 
$m$-dimensional vector uniformly distributed on the unit sphere 
$$S_1^m = \left\{ \bold x\in \Bbb R^m: \; \Vert\bold x\Vert^2 = x_1^2+ \dots +x_m^2=1 \right\} .$$
The motion is governed by a homogeneous Poisson process of rate $\lambda>0$ as follows. At each Poissonian instant, the particle instantaneously takes on a new random direction distributed uniformly on $S_1^m$ and keeps moving with the same speed $c$ 
until the next Poisson event occurs, then it takes on a new random direction again and so on. 

The particle's position $\bold X(t)=(X_1(t), \dots , X_m(t))$ at time $t>0$ is referred to as the $m$-dimensional {\it symmetric Markov random flight}. At arbitrary time instant $t>0$ the particle, with probability 1, is located in the closed $m$-dimensional ball of radius $ct$ centred at the origin $\bold 0$:
$$\mathcal B_{ct}^m = \left\{ \bold x\in \Bbb R^m : \; \Vert\bold x\Vert^2 = x_1^2+\dots+x_m^2\le c^2t^2 \right\} .$$

Denote by $\Phi(\bold x,t) = \text{Pr} \left\{ \bold X(t)\in d\bold x \right\}, \; \bold x\in\mathcal B_{ct}^m, \; t\ge 0,$  the distribution function of the process $\bold X(t)$, where $d\bold x\subset\Bbb R^m$ is an infinitesimal element in $\Bbb R^m$. At  arbitrary time instant $t>0$, the distribution function $\Phi(\bold x, t)$ consists of two components. 

The singular component is referred to the case when no Poisson events occur on the time interval $(0,t)$ and it is concentrated on the sphere 
$$S_{ct}^m =\partial\mathcal B_{ct}^m = \left\{ \bold x\in \Bbb R^m: \; \Vert\bold x\Vert^2 = x_1^2+\dots + x_m^2=c^2t^2 \right\} .$$
In this case, at time $t$, the particle is located on the sphere $S_{ct}^m$ and the probability of this event is: 
$\text{Pr} \left\{ \bold X(t)\in S_{ct}^m \right\} = e^{-\lambda t}$.

If at least one Poisson event occurs before time $t$, then the particle is located in the interior $\text{int} \; \mathcal B_{ct}^m$ 
of the ball $\mathcal B_{ct}^m$ and the probability of this event is: $\text{Pr} \left\{ \bold X(t)\in \text{int} \; \mathcal B_{ct}^m \right\} = 1 - e^{-\lambda t}$.
The part of the distribution function $\Phi(\bold x, t)$ corresponding to this case is concentrated strictly in the interior:  
$$\text{int} \; \mathcal B_{ct}^m = \left\{ \bold x\in \Bbb R^m: \; \Vert\bold x\Vert^2 = x_1^2 + \dots + x_m^2<c^2t^2 \right\}$$
of the ball $\mathcal B_{ct}^m$ and forms its absolutely continuous component. 

Let $p(\bold x,t), \; \bold x\in\mathcal B_{ct}^m , \; t>0,$ denote the probability density of distribution $\Phi(\bold x,t)$. 
Structurally, it has the form:  
$p(\bold x,t) = p^{(s)}(\bold x,t) + p^{(ac)}(\bold x,t) , \; \bold x\in\mathcal B_{ct}^m, \; t>0$,
where $p^{(s)}(\bold x,t)$ is the density (in the sense of generalized functions) of the singular component of $\Phi(\bold x, t)$ concentrated on $S_{ct}^m$ and $p^{(ac)}(\bold x,t)$ is the density of the absolutely continuous component of $\Phi(\bold x, t)$ concentrated in $\text{int} \; \mathcal B_{ct}^m$. 

The singular part of density $p(\bold x,t)$ is given by the formula: 
$$p^{(s)}(\bold x,t) =  \frac{e^{-\lambda t} \; \Gamma\left( \frac{m}{2} \right)}{2\pi^{m/2} (ct)^{m-1}} \;   
\delta(c^2t^2-\Vert\bold x\Vert^2) , \qquad \bold x\in\Bbb R^m, \;  t>0 ,$$
where $\delta(x)$ is the Dirac delta-function. This is the density of the uniform distribution on the surface of sphere $S_{ct}^m$.

The absolutely continuous part of density $p(\bold x,t)$ has the form: 
$$p^{(ac)}(\bold x,t) = f^{(ac)}(\bold x,t) \Theta(ct-\Vert\bold x\Vert) , \qquad \bold x\in\Bbb R^m, \; t>0,$$
where $f^{(ac)}(\bold x,t)$ is some positive function absolutely continuous in $\text{int} \; \mathcal B_{ct}^m$ and 
$\Theta(x)$ is the Heaviside unit-step function. 

Markov random flights in the Euclidean spaces of different dimensions were studied in a series of works \cite{frances,garcia,hughes,kol1,kol2,kol3,kol4,kol5,kol6,kolors,kolpin,mas,pin,sta1,sta2}. In the spaces $\Bbb R^2,\Bbb R^4,\Bbb R^6$ 
their distributions were, surprisingly, obtained in explicit forms \cite{kol2,kol5,kol6,kolors,mas,sta1,sta2}. These results will be used in the next section for deriving the stationary distributions of slow diffusion processes. However, first of all, we should determine a condition (in fact, a set of conditions) under which such stationary distributions exist.

Thus, we are interested in the conditions under which the $m$-dimensional Markov random flight $\bold X(t)$ can serve as an appropriate  mathematical model for a slow diffusion process $\bold D(t)$ in the Euclidean space $\Bbb R^m, \; m\ge 2$. First, since $\bold D(t)$ has a slow speed of propagation, then the random flight $\bold X(t)$ must have a slow speed too. Therefore, we should assume that the speed $c$ of process $\bold X(t)$ tends to zero. From this fact, as well as by taking into account some obvious physical reasonings, we can conclude that the intensity of switchings $\lambda$ should also tend to zero. It is clear that, in order to obtain important characteristics of the slow diffusion process $\bold D(t)$, we should consider random flight $\bold X(t)$ for large time $t$, that is for $t\to\infty$. This latter condition implies that we are interested in the stationary distribution of $\bold D(t)$. Since the parameters $\lambda$ and $c$ of the Markov random flight $\bold X(t)$ are connected with each other through the time, namely $\lambda$ 
is the mean number of changes of direction {\it per unit of time} and $c$ is the distance passed {\it per unit of time}, then we should assume that the products $\lambda t$ and $ct$ must tend to some finite limits. Note that $\lambda t$ is the mean number of changes of direction that have occured until time instant $t$, while $ct$ is the distance passed by time instant $t$. 

All these reasonings lead us to the following {\it slow diffusion condition} (SDC): 
\begin{equation}\label{SDC}
\aligned
& \lambda\to 0, \qquad c\to 0, \qquad t\to\infty , \\
&\quad \lambda t \to a > 0 , \qquad ct\to\varrho > 0. 
\endaligned 
\end{equation}

Condition (\ref{SDC}) has very clear physical sense. A slow diffusion process $\bold D(t)$ can be simulated by the symmetric 
Markov random flight $\bold X(t)$ with small speed and intensity of switchings. This smallness of the values of the parameters $c$ and $\lambda$ of $\bold X(t)$ implies that we should study the process on very long time intervals on which its probabilistic characteristics become stable and tend to the stationary ones. In the next section we will show that SDC (\ref{SDC}) leads to the stationary distributions indeed and, therefore, such random flights can be applied for modelling the slow diffusion processes in the Euclidean spaces of arbitrary dimension.

\section{Stationary Densities in Low Dimensions}
\numberwithin{equation}{section} 

Let, as above, $p(\bold x,t), \; \bold x\in\Bbb R^m, \;  t>0,$ denote the transition density of the symmetric Markov 
random flight $\bold X(t)$. In this section we derive, under the SDC (\ref{SDC}), the densities  
$$q(\bold x) = \lim_{\substack{c,\lambda\to 0, \; t\to\infty\\\lambda t\to a, \; ct\to\varrho}} p(\bold x,t)$$
of the stationary distributions of $\bold X(t)$ in some Euclidean spaces of low dimensions that 
can be treated as the stationary densities of the slow diffusion processes. The derivation is based on the closed-form expressions 
for the transition densities of $\bold X(t)$ obtained in some previous works. In the space $\Bbb R^3$ the stationary density 
is presented in an asymptotic form. 

\vskip 0.2cm

{\bf Stationary Density on the Line.}
The one-dimensional Markov random flight is represented by the Goldstein-Kac telegraph process described above. Its density is given by formula (\ref{dens1}). Taking into account the well known relation $I_0'(z)=I_1(z)$ connecting the modified Bessel functions of the zero and first orders, we can rewrite density (\ref{dens1}) in the alternative form (see also \cite[formula (2.5.3)]{kolrat}): 
\begin{equation}\label{densalter}
\aligned p(x,t) & = \frac{e^{-\lambda t}}{2} \left[ \delta(ct+x) +
\delta(ct-x) \right]\\
& + \frac{\lambda e^{-\lambda t}}{2c} \left[ I_0\left(
\frac{\lambda}{c} \sqrt{c^2t^2-x^2}\right) + \frac{c
t}{\sqrt{c^2t^2-x^2}} \; I_1\left( \frac{\lambda}{c}
\sqrt{c^2t^2-x^2}\right)
\right] \Theta(ct-\vert x\vert) , \\
& \qquad\qquad\qquad x\in (-\infty, \infty), \qquad t>0, 
\endaligned
\end{equation}
where
$$I_1(z) = \sum_{k=0}^{\infty} \frac{1}{k! \; (k+1)!} \left( \frac{z}{2} \right)^{2k+1}$$
is the modified Bessel function of first order and $\Theta(x)$ is the Heaviside unit-step function given by (\ref{heaviside}). 

Under the SDC (\ref{SDC}), density (\ref{densalter}) transforms into the stationary density: 
\begin{equation}\label{statdens1}
\aligned q(x) & = \frac{e^{-a}}{2} \left[ \delta(\varrho+x) +
\delta(\varrho-x) \right]\\
& + \frac{a e^{-a}}{2\varrho} \left[ I_0\left(
\frac{a}{\varrho} \sqrt{\varrho^2-x^2}\right) + \frac{\varrho}{\sqrt{\varrho^2-x^2}} \; 
I_1\left( \frac{a}{\varrho} \sqrt{\varrho^2-x^2}\right) \right] \Theta(\varrho-\vert x\vert) , \\
& \qquad\qquad\qquad x\in (-\infty, \infty), \quad a>0, \quad \varrho >0. 
\endaligned
\end{equation} 

The shape of the absolutely continuous part of stationary density $q(x)$ given by formula (\ref{statdens1}) is presented in Fig. 1. 
\begin{center}
\begin{figure}[htbp]
\centerline{\includegraphics[width=12cm,height=7cm]{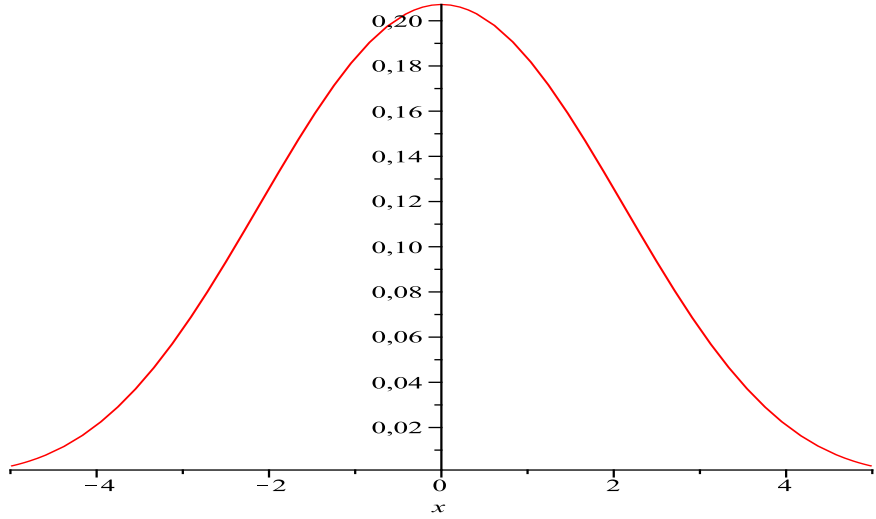}}
\caption{\it The shape of stationary density (\ref{statdens1}) (for $\varrho=5, \; a =7, \; |x|<5$)}
\end{figure}
\end{center}
One can see that, for increasing $a$ and fixed $\varrho$, this curve becomes more and more peaked. On the other hand, for increasing $\varrho$ and fixed $a$, this graphics becomes more even with decreasing peak.

\vskip 0.2cm

{\bf Stationary Density in the Plane.}
The transition density of the symmetric Markov random flight $\bold X(t)$ in the Euclidean plane $\Bbb R^2$ 
is given by the formula (see \cite{kolors,mas,sta1,sta2}):
\begin{equation}\label{dens2}
p(\bold x,t) = \frac{e^{-\lambda t}}{2\pi ct} \; \delta(c^2t^2-\Vert\bold x\Vert^2) + \frac{\lambda}{2\pi c} \; \frac{\exp \left( -\lambda t + \frac{\lambda}{c} \sqrt{c^2t^2-\Vert\bold x\Vert^2} \right)}{\sqrt{c^2t^2-\Vert\bold x\Vert^2}} \;
\Theta(ct-\Vert\bold x\Vert), 
\end{equation}
$$\bold x=(x_1,x_2)\in\Bbb R^2 , \quad \Vert\bold x\Vert = \sqrt{x_1^2+x_2^2} , \quad t>0 .$$

Under the SDC (\ref{SDC}), density (\ref{dens2}) transforms into the stationary density: 
\begin{equation}\label{statdens2}
q(\bold x) = \frac{e^{-a}}{2\pi \varrho} \; \delta(\varrho^2-\Vert\bold x\Vert^2) + \frac{a}{2\pi \varrho} \; \frac{\exp \left( -a + \frac{a}{\varrho} \sqrt{\varrho^2-\Vert\bold x\Vert^2} \right)}{\sqrt{\varrho^2-\Vert\bold x\Vert^2}} \;
\Theta(\varrho-\Vert\bold x\Vert) ,   
\end{equation}
$$\bold x\in\Bbb R^2, \quad a>0, \quad \varrho >0 .$$

\begin{center}
\begin{figure}[htbp]
\centerline{\includegraphics[width=12cm,height=7cm]{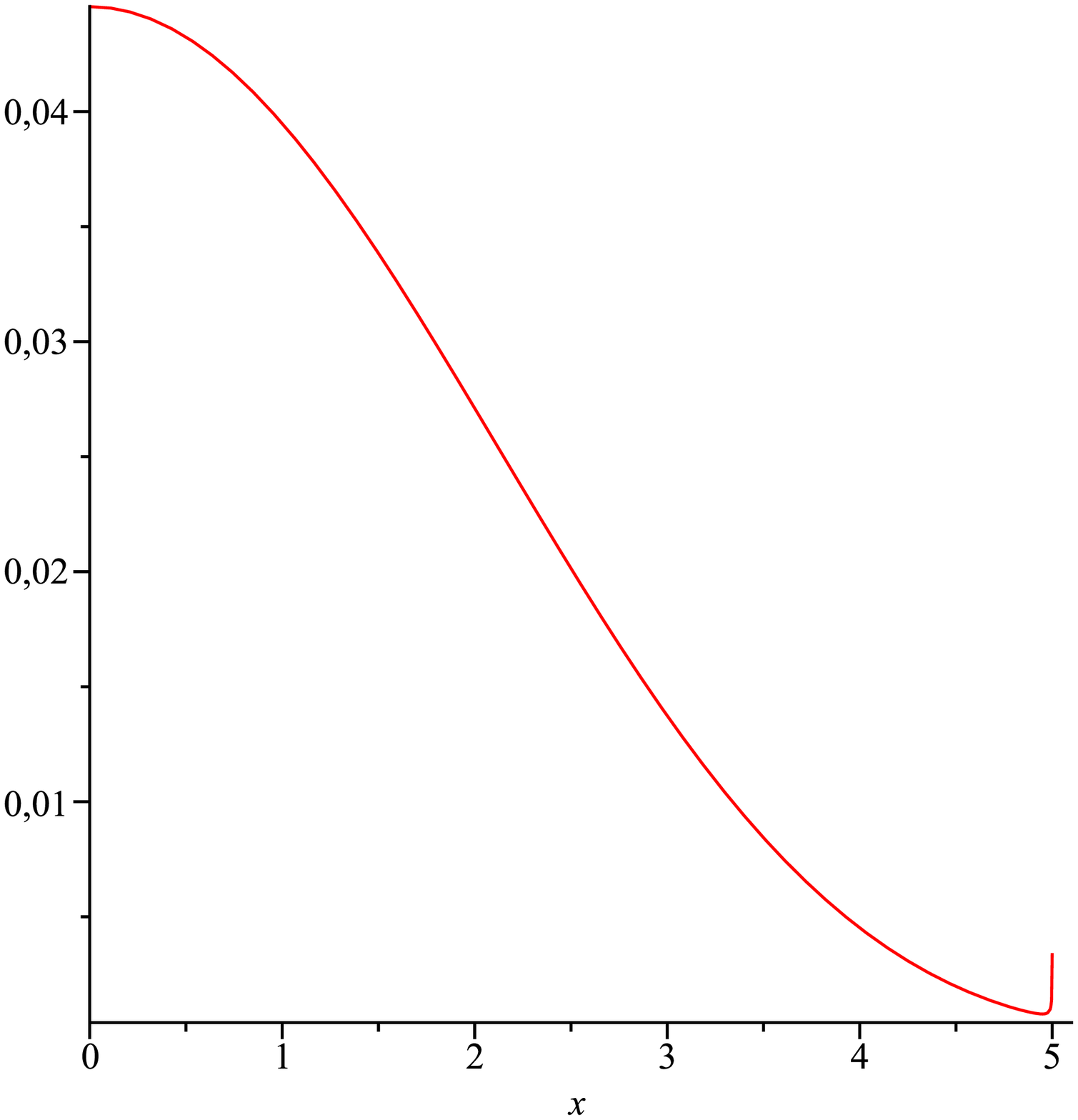}}
\caption{\it The shape of the section of stationary density (\ref{statdens2}) (for $\varrho=5, \; a =7, \; \Vert\bold x\Vert<5$)}
\end{figure}
\end{center} 

The shape of the section of the absolutely continuous part of stationary density $q(\bold x)$ given by formula (\ref{statdens2}) 
is plotted in Fig. 2.
It shows the behaviour of $q(\bold x)$ as the distance $\Vert\bold x\Vert$ from the origin $\bold 0\in\Bbb R^2$ grows. We see that 
$q(\bold x)$ has local maximum at the origin $\bold 0$ and it is decreasing, as $\Vert\bold x\Vert$ grows. Note also that density $q(\bold x)$ becomes infinite near the border $\Vert\bold x\Vert=\varrho$, that is, $\lim\limits_{\Vert\bold x\Vert\to\varrho -0} q(\bold x) = \infty$. This follows from the form of stationary density (\ref{statdens2}).

\vskip 0.2cm

{\bf Stationary Density in the Space $\Bbb R^4$.} 
The transition density of the symmetric Markov random flight $\bold X(t)$ in the four-dimensional space $\Bbb R^4$ 
has the form (see \cite[formula (19)]{kol6}):
\begin{equation}\label{dens4} 
\aligned 
p(\bold x,t) & = \frac{e^{-\lambda t}}{2\pi^2 (ct)^3} \; \delta(c^2t^2-\Vert\bold x\Vert^2) \\ 
& + \frac{\lambda t}{\pi^2 (ct)^4} \left[ 2 + \lambda t \left( 1 -
\frac{\Vert\bold x \Vert^2}{c^2t^2} \right) \right] \exp \left( - \frac{\lambda}{c^2 t} \; \Vert\bold x \Vert^2 \right) \; 
\Theta(ct-\Vert\bold x\Vert) , \endaligned 
\end{equation}
$$\bold x = (x_1,x_2,x_3,x_4)\in\Bbb R^4, \quad \Vert\bold x \Vert = \sqrt{x_1^2+x_2^2+x_3^2+x_4^2}, \quad t>0.$$

Under the SDC (\ref{SDC}), density (\ref{dens4}) transforms into the stationary density: 
\begin{equation}\label{statdens4}
\aligned 
q(\bold x) & = \frac{e^{-a}}{2\pi^2 \varrho^3} \; \delta(\varrho^2-\Vert\bold x\Vert^2) \\ 
& + \frac{a}{\pi^2 \varrho^4} \left[ 2 + a \left( 1 -
\frac{\Vert\bold x \Vert^2}{\varrho^2} \right) \right] \exp \left( - \frac{a}{\varrho^2} \; \Vert\bold x \Vert^2 \right) \; 
\Theta(\varrho-\Vert\bold x\Vert) , \endaligned 
\end{equation}
$$\bold x\in\Bbb R^4, \quad a>0, \quad \varrho >0 .$$ 

The shape of the section of the absolutely continuous part of stationary density $q(\bold x)$ given by formula (\ref{statdens4}) 
is presented in Fig. 3.

\begin{center}
\begin{figure}[htbp]
\centerline{\includegraphics[width=12cm,height=7cm]{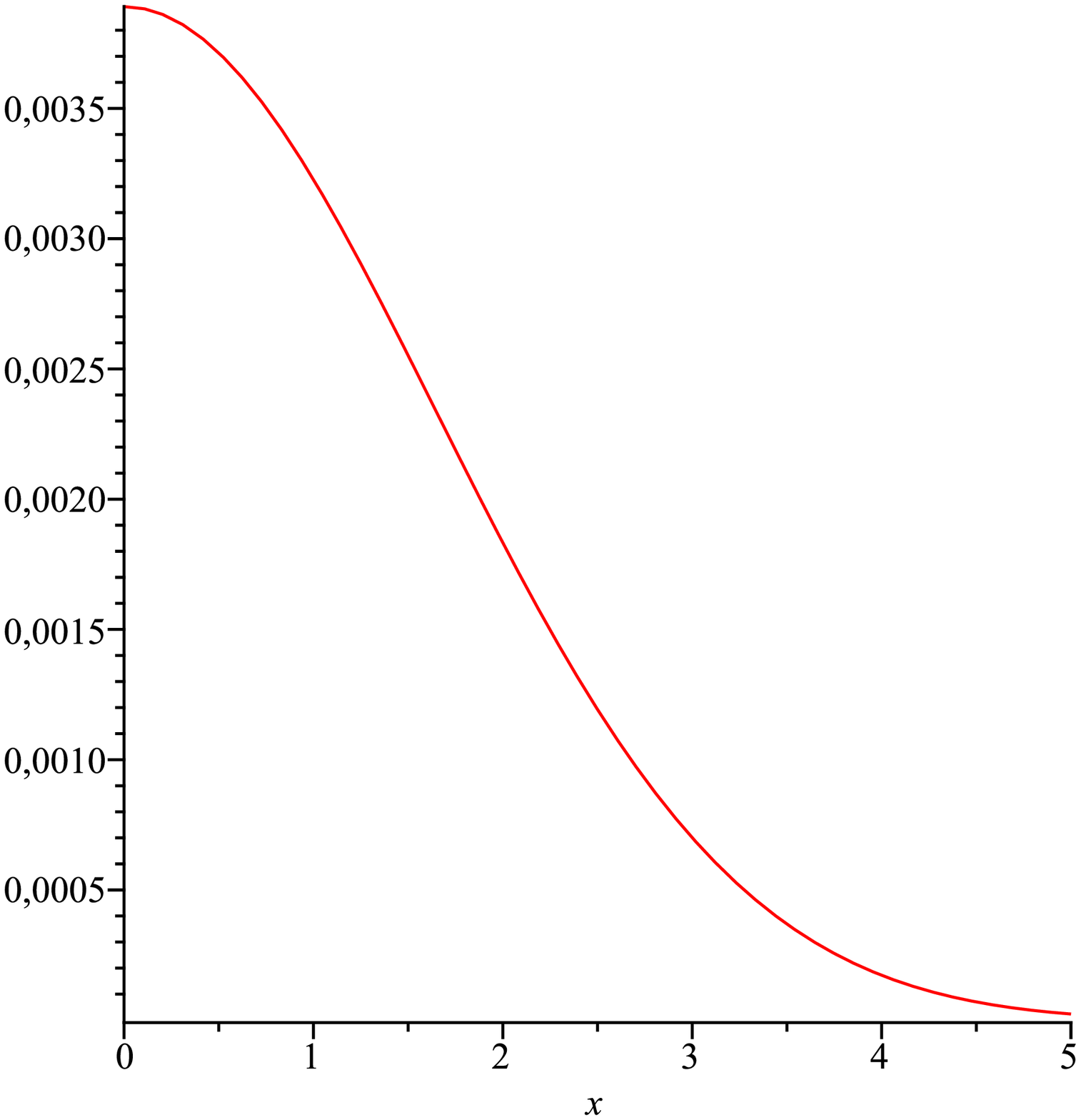}}
\caption{\it The shape of the section of stationary density (\ref{statdens4}) (for $\varrho=5, \; a =4, \; \Vert\bold x\Vert<5$)}
\end{figure}
\end{center}

We see that stationary density (\ref{statdens4}) takes maximal value at the origin and decreases smoothly up to the boundary. It is continuous and takes its minimal value on the boundary of the diffusion area.  

\vskip 0.2cm

{\bf Stationary Density in the Space $\Bbb R^6$.} 
The transition density of the symmetric Markov random flight $\bold X(t)$ in the six-dimensional space $\Bbb R^6$ 
has the form (see \cite[the Theorem]{kol2}):
\begin{equation}\label{dens6} 
\aligned 
p(\bold x,t) & = \frac{e^{-\lambda t}}{\pi^3 (ct)^5} \; \delta(c^2t^2-\Vert\bold x\Vert^2) + 
\biggl[ \frac{16\lambda t \; e^{-\lambda t}}{\pi^3 (ct)^6} \left( 1 - \frac{5}{6} \;
\frac{\Vert \bold x \Vert^2}{c^2t^2} \right) \\
& + \frac{e^{-\lambda t}}{2\pi^3 (ct)^6} \sum_{n=2}^{\infty}
(\lambda t)^n (n+1)! \; \sum_{k=0}^{n+1}
\frac{(k+1)(k+2)(n+2k+1)}{3^k (n-k+1)! (n+k-2)!} \\
& \qquad\qquad\qquad\qquad \times F \left( -(n+k-2), k+3; \; 3; \frac{\Vert \bold x \Vert^2}{c^2t^2} \right) \biggr]
\Theta(ct-\Vert\bold x\Vert) ,
\endaligned 
\end{equation}
$$\bold x = (x_1,x_2,x_3,x_4,x_5,x_6)\in\Bbb R^6, \quad \Vert\bold x \Vert = \sqrt{x_1^2+x_2^2+x_3^2+x_4^2+x_5^2+x_6^2}, \quad t>0,$$
where 
$$F(\alpha,\beta; \; \gamma; \; z) = \sum_{k=0}^{\infty} \frac{(\alpha)_k \; (\beta)_k}{(\gamma)_k} \; \frac{z^k}{k!}$$ 
is the Gauss hypergeometric function. 

Under the SDC (\ref{SDC}), density (\ref{dens4}) transforms into the stationary density: 
\begin{equation}\label{statdens6}
\aligned 
q(\bold x) & = \frac{e^{-a}}{\pi^3 \varrho^5} \; \delta(\varrho^2-\Vert\bold x\Vert^2) + 
\biggl[ \frac{16 a \; e^{-a}}{\pi^3 \varrho^6} \left( 1 - \frac{5}{6} \; \frac{\Vert \bold x \Vert^2}{\varrho^2} \right) \\ 
& + \frac{e^{-a}}{2\pi^3 \varrho^6} \sum_{n=2}^{\infty} a^n (n+1)! \; 
\sum_{k=0}^{n+1} \frac{(k+1)(k+2)(n+2k+1)}{3^k (n-k+1)! (n+k-2)!} \\
& \qquad\qquad\qquad\qquad \times F \left( -(n+k-2), k+3; \; 3; \frac{\Vert \bold x \Vert^2}{\varrho^2} \right) \biggr]
\Theta(\varrho-\Vert\bold x\Vert) ,
\endaligned 
\end{equation}
$$\bold x\in\Bbb R^6, \quad a>0, \quad \varrho >0 .$$ 

The shape of the section of the absolutely continuous part of stationary density $q(\bold x)$ given by formula (\ref{statdens6}) 
is plotted in Fig. 4.

\begin{center}
\begin{figure}[htbp]
\centerline{\includegraphics[width=12cm,height=7cm]{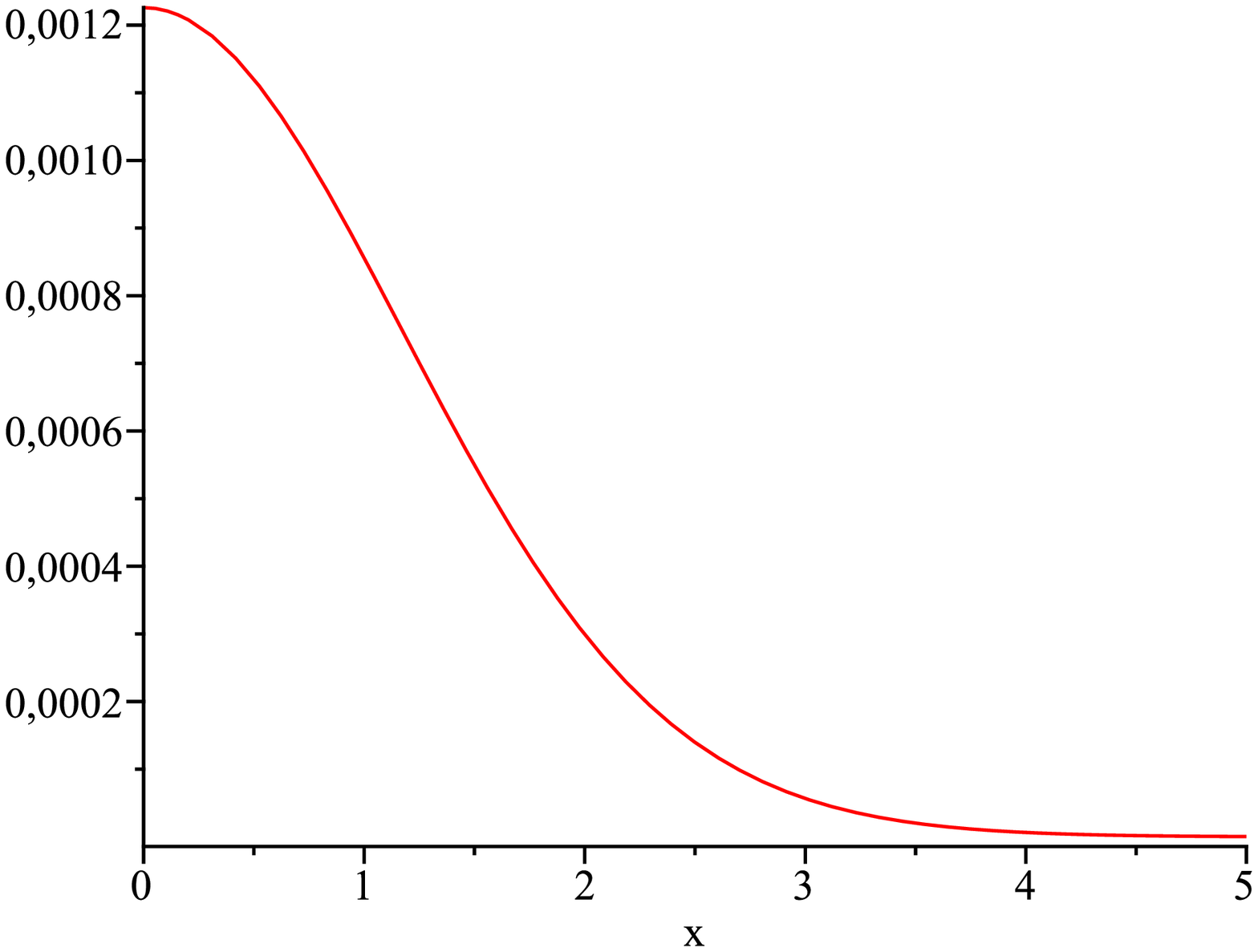}}
\caption{\it The shape of the section of stationary density (\ref{statdens6}) (for $\varrho=5, \; a =4, \; \Vert\bold x\Vert<5$)}
\end{figure}
\end{center} 

This figure shows that the six-dimensional stationary density (\ref{statdens6}) behaves very likely to its four-dimensional counterpart 
(see Fig. 3), namely, it takes maximal value at the origin and decreases smoothly up to the boundary. On the boundary it is continuous and takes its minimal value. As we will see in the next subsection, in the three-dimensional space $\Bbb R^3$, under an additional condition of smallness of the parameter $a$, the behaviour of the stationary density changes drastically. 

\vskip 0.2cm

{\bf Asymptotic Stationary Density in the Space $\Bbb R^3$.}
The above stationary densities of the Markov random flights in the spaces $\Bbb R^1, \Bbb R^2, \Bbb R^4, \Bbb R^6$ were derived from 
the explicit forms of their transition densities in these spaces obtained in a series of previous works \cite{kol2,kol5,kol6,kolors,kolrat,mas,sta1,sta2}. As to the finite-velocity random motions driven by a Poisson process in the three-dimensional space $\Bbb R^3$ is concerned, the situation is more complicated. This motion was the subject of research in some works \cite{kol3,mas,sta1} and a few interesting results were obtained. In particular, the explicit form of the Laplace-Fourier  transform of the transition density of the three-dimensional symmetric Markov random flight was derived by different methods 
(see \cite[p. 1054]{kol3}, \cite[formula (45)]{mas}, \cite[formula (5.8)]{sta1}). However, the problem of inverting this 
Laplace-Fourier transform and obtaining a closed-form expression for the density remains still open. That is why an alternative 
approach was developed in \cite{kol1} yielding an asymptotic formula for the transition density of the symmetric Markov random flight 
in the space $\Bbb R^3$ that will be used in this subsection for obtaining an asymptotic stationary density of the process. 

Let $\bold X(t) = (X_1(t), X_2(t), X_3(t)), \; t>0$ be the symmetric Markov random flight in the Euclidean space $\Bbb R^3$ with constant speed $c>0$ and the intensity of switchings $\lambda>0$. Let $p(\bold x,t), \; \bold x\in\Bbb R^3, \;  t>0,$ denote the transition density of $\bold X(t)$. It was shown in \cite[Theorem 3 and comments on page 447]{kol1} that, for arbitrary $t>0$, the following asymptotic relation holds:  
\begin{equation}\label{dens3}
\aligned 
p(\bold x,t) & = \frac{e^{-\lambda t}}{4\pi (ct)^2} \; \delta(c^2t^2-\Vert\bold x\Vert^2) + e^{-\lambda t} \biggl[ \frac{\lambda}{4\pi c^2 t \; \Vert\bold x\Vert} \ln\left( \frac{ct+\Vert\bold x\Vert}{ct-\Vert\bold x\Vert} \right) \\ 
& \qquad\qquad\qquad\qquad  + \frac{\lambda^2}{2\pi^2 c^2 \; \sqrt{c^2t^2-\Vert\bold x\Vert^2}} + \frac{\lambda^3}{8\pi c^3} \biggr] \Theta(ct-\Vert\bold x\Vert) + o((\lambda t)^3) , 
\endaligned
\end{equation}
$$\bold x = (x_1, x_2, x_3) \in\Bbb R^3, \qquad \Vert\bold x\Vert = \sqrt{x_1^2+x_2^2+x_3^2}, \qquad t>0.$$

Therefore, under the SDC (\ref{SDC}) and the additional condition $0<a\ll 1$, density (\ref{dens3}) transforms into the asymptotic 
stationary density: 
\begin{equation}\label{statdens3} 
\aligned 
q(\bold x) & = \frac{e^{-a}}{4\pi \varrho^2} \; \delta(\varrho^2-\Vert\bold x\Vert^2) + e^{-a} \biggl[ \frac{a}{4\pi \varrho^2 \;  \Vert\bold x\Vert} \ln\left( \frac{\varrho+\Vert\bold x\Vert}{\varrho-\Vert\bold x\Vert} \right) \\ 
& \qquad\qquad\qquad\qquad  + \frac{a^2}{2\pi^2 \varrho^2 \; \sqrt{\varrho^2-\Vert\bold x\Vert^2}} + \frac{a^3}{8\pi \varrho^3} \biggr] \Theta(\varrho-\Vert\bold x\Vert) + o(a^3) ,
\endaligned
\end{equation}
$$\bold x\in\Bbb R^3, \quad 0<a\ll 1, \quad \varrho >0 .$$

The shape of the section of the absolutely continuous part of asymptotic stationary density $q(\bold x)$ given by formula (\ref{statdens3}) (for $\varrho=5, \; a =0.01, \; \Vert\bold x\Vert<5$), is presented in Fig. 5. The error of calculations in this graphics does not exceed $10^{-6}$. 

\begin{center}
\begin{figure}[htbp]
\centerline{\includegraphics[width=12cm,height=7cm]{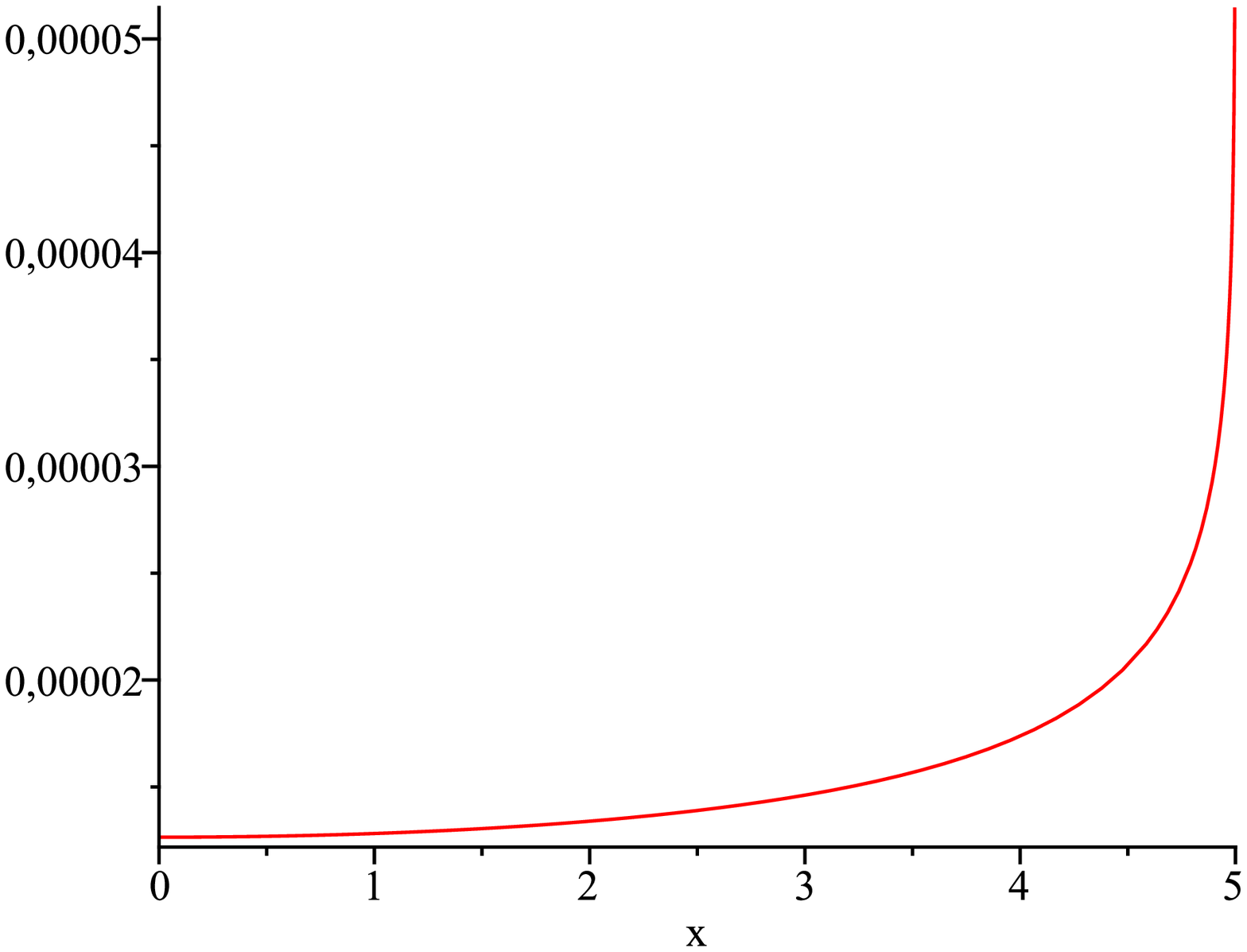}}
\caption{\it The shape of the section of stationary density (\ref{statdens3}) (for $\varrho=5, \; a =0.01, \; \Vert\bold x\Vert<5$)}
\end{figure}
\end{center} 

We see the essential difference in the behaviour of the three-dimensional density and those in the even-dimensional spaces $\Bbb R^2, \Bbb R^4, \Bbb R^6$. All the densities in these even-dimensional spaces have almost the same behaviour, namely, they are mostly concentrated near the origin and are decreasing, as the distance from the origin grows. Near the border $\Vert\bold x\Vert=\varrho$ these densities take minimal values (except the two-dimensional density that is infinite on the border, but the probability of being near the border is, nevertheless, small like in other even-dimensional cases). In contrast to these even-dimensional spaces, the stationary density in the space $\Bbb R^3$ takes {\it minimal} value at the origin and is slowly growing, as the distance from the origin grows. Thus, one can say that there is a rarefaction area near the origin. When approaching the border, the density begins to rise sharply and it becomes infinite on the border. Therefore, one can conclude that the overwhelming part of density is concentrated near the border and this fact looks very similar to the behaviour of the fundamental solution (the Green's function) of the three-dimensional wave equation which is concentrated on the surface of the diffusion sphere (see, for instance, \cite[section 11, subsection 7]{vlad}). This analogy becomes particularly relevant if we take into account that random flights are, in fact, a sort of wave processes and, therefore, one can expect that they should possess some their main properties. In particular, the well-known difference in the behaviour of wave processes in the Euclidean spaces of even and odd dimensions caused by the feasibility of the Huygens principle can, apparently, be a reasonable explanation of the difference in the behaviour of the densities of Markov random flights in the spaces of even and odd dimensions.

\bigskip 

\begin{center}
\bf Conclusions 
\end{center}

The presented conception of slow diffusion processes based on Markov random flights gives a tool for describing their 
distributions on long time intervals in the Euclidean spaces of low, most important, dimensions. The core of this conception 
is the slow diffusion condition (\ref{SDC}) connecting the speed of propagation and the intensity of switchings of the process with 
the time interval of its evolution. On long time intervals, such distributions tend to the respective stationary distributions and this is valid in any dimension (asymptotically in the space $\Bbb R^3$). The slow diffusion condition (\ref{SDC}) can, therefore, be considered as the slow-velocity counterpart of the classical Kac's fast diffusion condition (\ref{KacCond}).

\end{document}